\documentclass[letter]{article}
\usepackage{amsmath}
\usepackage{amsfonts}
\usepackage{amssymb}
\usepackage{graphicx}
\usepackage{enumerate}
\usepackage{url}

\begin{document}

\renewcommand{\thefootnote}{\fnsymbol{footnote}}

\newtheorem{theorem}{Theorem}[section]
\newtheorem{corollary}[theorem]{Corollary}
\newtheorem{definition}[theorem]{Definition}
\newtheorem{conjecture}[theorem]{Conjecture}
\newtheorem{question}[theorem]{Question}
\newtheorem{lemma}[theorem]{Lemma}
\newtheorem{proposition}[theorem]{Proposition}
\newtheorem{example}[theorem]{Example}
\newenvironment{proof}{\noindent {\bf
Proof.}}{\rule{3mm}{3mm}\par\medskip}
\newcommand{\remark}{\medskip\par\noindent {\bf Remark.~~}}
\newcommand{\pp}{{\it p.}}
\newcommand{\de}{\em}

\title{ \bf The largest $n-1$ Hosoya indices of unicyclic graphs
\thanks{
This paper was supported by Foundation of Education Department of Shandong Province (J07YH03),
NNSFC (10871205, 70901048), NSFSD (No. Y2008A04, BS2010SF017) and Research projects 174010 and
174033 of the Serbian Ministry of Science.}}

\author{Guihai Yu, \ Lihua Feng\\
{\small \it School of Mathematics, Shandong Institute of Business and Technology} \\
{\small \it 191 Binhaizhong Road, Yantai, Shandong, P.R. China, 264005.}\\
{\small e-mail: { \tt yuguihai@126.com, fenglh@163.com}} \\
\and
Aleksandar Ili\' c \footnotemark[3] \\
{\small \it Faculty of Sciences and Mathematics} \\
{\small \it University of Ni\v s, Vi\v segradska 33, 18000 Ni\v s, Serbia} \\
{\small e-mail: { \tt aleksandari@gmail.com}}\\
}

\date{October 1, 2010}

\maketitle
\vspace{-0.5cm}

\begin{abstract}
The Hosoya index $Z (G)$ of a graph $G$ is defined as the total number of edge independent sets of
$G$. In this paper, we extend the research of [J. Ou,  On extremal unicyclic molecular graphs with
maximal Hosoya index, \textit{Discrete Appl. Math.} 157 (2009) 391--397.] and [Y. Ye, X. Pan, H.
Liu, Ordering unicyclic graphs with respect to Hosoya indices and Merrifield-Simmons indices,
\textit{MATCH Commun. Math. Comput. Chem.} 59 (2008) 191--202.] and order the largest $n-1$
unicyclic graphs with respect to the Hosoya index.
\end{abstract}

\vskip 0.1cm

\section{Introduction}

In this paper, we follow the standard notation in graph theory in \cite{Bondy}. Let $G = (V, E)$ be
a simple connected graph of order $n$. Two distinct edges in a graph $G$ are {\it independent} if
they are not incident with a common vertex in $G$. A set of pairwise independent edges in $G$ is
called a {\it matching}. A $k-$matching of $G$ is a set of $k$ mutually independent edges.

In theoretical chemistry molecular structure descriptors are used for modeling physico-chemical,
pharmacologic, toxicologic, biological and other properties of chemical compounds. The {\it Hosoya
index} $Z (G)$ of a graph $G$ is defined as the total number of its matchings {\cite{Hosoya}}. If
$m(G,k)$ denotes the number of its $k-$matchings, then
$$Z(G)=\sum_{k=0}^{\lfloor n / 2 \rfloor}m(G,k).$$

The Hosoya index has been studied intensively in the literature (see survey \cite{WaGu10} and
references \cite{Chan,Gutman01,hua,li09,Wagner,xu1,zhang04}). It is shown in \cite{Gutman 03} that the linear
hexagonal chain is the unique graph with minimal Hosoya index among all hexagonal chains, while in
\cite{zhang01} and \cite{zhang02} the authors proved that zig-zag hexagonal chain is the unique
chain with the maximal Hosoya index among all hexagonal chains. As for trees, it has been shown
that the path has the maximal Hosoya index and the star has the minimal Hosoya index \cite{Gutman}.
Recently, Hou \cite{Hou} characterized the trees with a given size of matching having minimal and
second minimal Hosoya index. Pan et al. in \cite{pan} characterized the trees with given diameter
and having minimal Hosoya index. In \cite{Ou} Ou used linear algebra theory on permanents to
characterize the minimal and second minimal Hosoya index of unicyclic graphs with fixed girth. In
\cite{xu} the authors considered the Hosoya and Merrifield-Simmons indices of unicyclic graphs with
maximum degree $\Delta$ and characterized the graphs with the maximal Hosoya index and the minimal
Merrifield-Simmons index. \vspace{0.1cm}

For $v \in V (G)$, we denote by $G - v$ the graph obtained from $G$ by deleting the vertex $v$
together with their incident edges. For $e \in E (G)$, we denote by $G - e$ the graph obtained from
$G$ by removing the edge $e$. Let $deg (v)$ denotes the vertex degree of $v$. We denote by $P_n$,
$S_n$ and $C_n$ the path, the star and the cycle on $n$ vertices, respectively. By $L_{n,k}$ we
denote the graph obtained from $C_{k}$ and $P_{n-k+1}$ by identifying a vertex of $C_{k}$ with one
end vertex of $P_{n-k+1}$.

Let $F_{n}$ denotes the $n$-th Fibonacci number. It is well-known that $F_{n}$ satisfy the
following recursive relations:
\begin{equation}
F_{n}=F_{n-1}+F_{n-2}, \qquad  F_{1}=1,\ F_{0}=0, \ n \geq 2
\end{equation}
\begin{equation}
F_{n}=F_{k}F_{n-k+1}+F_{k-1}F_{n-k},\qquad   1\leq k\leq n
\end{equation}
\begin{equation}
\label{fi3} F_{m}F_{n+1}-F_{n}F_{m+1}=(-1)^n F_{m-n}
\end{equation}
\begin{equation}
F_{-n}=(-1)^{n+1} F_n, \qquad  n\geq 0.
\end{equation}

Let ${\cal{U}}_{n,k}$ be the set of all unicyclic graphs of order $n\geq 3$ with girth $k\geq 3$.
It is easy to see that if $k=n$ or $n-1$, there is only one unicyclic graph. Therefore, we can
assume that $3\leq k\leq n-2$. Deng et al. in \cite{Deng} showed that the unique extremal unicyclic
graph in ${\cal{U}}_{n,k}$ with maximal Hosoya index is $L_{n,k}$. Ou in \cite{Ou jianping} proved
the following:
\begin{theorem}\label{ou}
Let $G$ be a connected unicyclic graph of order $n \geq 5$. If $G \not \cong C_n$, then $Z (G) \leq
F_{n + 1} + 2 F_{n - 3}$, with the equality holding if and only if $G \cong L_{n, 4}$ or $G \cong
L_{n, n - 2}$.
\end{theorem}

In \cite{ye}, the authors presented the smallest $\lfloor \frac{k}{2}\rfloor+2$ Hosoya indices of
unicyclic graphs with given girth $k$. A fully loaded unicyclic graph is a unicyclic graph with the
property that there is no vertex with degree less than 3 in its unique cycle. Hua in \cite{hua1}
determined the extremal fully loaded unicyclic graphs with minimal, second-minimal and
third-minimal Hosoya indices. Deng in \cite{Deng08} obtained the largest Hosoya index of $(n, n +
1)$ graphs. Inspired by these results, in this paper we determine the first $n-1$ unicyclic graphs
with the largest Hosoya index among all unicyclic graphs of order $n$.

This paper is organized as follows. In Section 2, we give basic results concerning Hosoya index and
some useful lemmas. In Section 3, we characterize unicyclic graphs in ${\cal{U}}_{n,k}$ with the
second maximal Hosoya index. Finally in Section 4, we order the largest $n-1$ unicyclic $n$-vertex
graphs with respect to the Hosoya index.

\section{Preliminary results}

We list some basic results from \cite{Gutman} concerning Hosoya index of a graph $G$:
\begin{enumerate}[($i$)]
\item  Let $e=uv$ be an edge of a graph $G$. Then
\begin{equation}
\label{formula4} Z(G)=Z(G-e)+Z(G-uv).
\end{equation}

\item Let $v$ be a vertex of a graph $G$. Then
\begin{equation}
\label{formula6} Z(G)=Z(G-v)+\sum_{u \sim v}Z(G-uv),
\end{equation}
where the summation extends over all vertices adjacent to $v$.

In particular, if $v$ is a pendant vertex of $G$ and $u$ is the only vertex adjacent to $v$, we
have $Z(G)=Z(G-v)+Z(G-u-v)$.

\item If $G_{1},G_{2},\ldots,G_{k}$ are connected components of $G$, then
$Z(G)=\prod_{i=1}^{k}Z(G_{i}).$

\item For paths, stars and cycles, we have
$$Z(P_{0})=0, \ Z(P_{1})=1, \qquad Z(P_{n})=F_{n+1}, \ n\geq2.$$
$$Z(S_{n})=n; \qquad Z(C_{n})=F_{n-1}+F_{n+1}.$$
\end{enumerate}

\begin{lemma} \label{1}\cite{Li}
Let $T_{n}$ be a tree of order $n$. Then $Z(S_{n})\leq Z(T_{n})\leq
Z(P_{n}),$ with the left equality holding if and only if $T_{n}\cong
S_{n}$ and the right equality holding if and only if $T_{n}\cong P_{n}$.
\end{lemma}

\begin{lemma} \label{4} \cite{Deng}
Let $U_{n,k}\in {\cal{U}}_{n,k}$ be a unicyclic graph with girth
$k\geq 3$. Then $Z(U_{n,k})\leq Z(L_{n,k}) = F_{n+1} + F_{k-1} F_{n-k+1},$ with the equality holding if
and only if $U_{n,k}\cong L_{n,k}$.
\end{lemma}

\begin{lemma}\cite{Wagner}\label{5}
Let $G$  be a connected graph and  $v\in V(G)$. Suppose $P(n,k,G,v)$ denotes the graph obtained
from $G$ by identifying $v$ with the vertex $v_{k}$ of a simple path $v_{1}v_{2}\cdots v_{n}$ (see
Fig.~1). Let $n=4m+i$ $(i\in \{0,1,2,3\}, m \geq 0)$. Then
$$Z(P(n,2,G,v))<Z(P(n,4,G,v))<\cdots<Z(P(n,2m,G,v))<$$
$$< Z(P(n,2m-1+2l,G,v))<\cdots<Z(P(n,3,G,v))<Z(P(n,1,G,v)),$$
where $l = \lfloor \frac{i}{2} \rfloor$.
\end{lemma}

\setlength{\unitlength}{0.05 in}
\begin{center}
\begin{picture}(80,20)
\put(40,13.5){\circle{20}} \put(40,8){\circle*{1}}
\put(40,8){\line(-10,0){5}}\put(40,8){\line(10,0){5}} \put(40,9){$v$}
\put(35,8){\circle*{1}}\put(25,8){$\ldots$}\put(30,8){$\ldots$}\put(25,8){\circle*{1}}\put(20,8){\circle*{1}}
\put(20,5){$v_{1}$}\put(20,8){\line(10,0){5}}
\put(40,5){$v_{k}$}\put(40,15){$G$}\put(45,8){\circle*{1}}\put(45,8){$\ldots$}
\put(50,8){$\ldots$}\put(55,8){\circle*{1}}\put(55,8){\line(10,0){5}} \put(60,8){\circle*{1}}
\put(60,5){$v_{n}$} \put(18,0){\textbf{Figure 1.} The graph $P (n, k, G, v)$. }
\end{picture}
\end{center}

Similarly, one can verify that for $n = 4m + i$ $(i \in \{0,1,2,3\})$ and $l = \lfloor \frac{i}{2}
\rfloor$, the following chain of inequalities holds
\begin{equation}
\label{fib} F_0 F_n < F_2 F_{n - 2} < F_4 F_{n - 4} <\cdots <F_{2m} F_{2m + i} < F_{2m - 1 + 2l}
F_{2m + i + 1 - 2l} < \cdots < F_{3} F_{n - 3} <F_1 F_{n - 1}.
 \end{equation}

By repeated use of Lemma \ref{5}, the extremal unicyclic graph with the maximal Hosoya index has
only paths attached to some vertices of a cycle.

\begin{lemma}\cite{Deng}\label{6}
Let $P=uu_1u_2\cdots u_tv$ be a path in $G$ and $G\neq P$, the degrees of $u_1, u_2, \ldots, u_t$
in $G$ are 2. Let $G_1$ denotes the graph that results from identifying $u$ with the vertex $v_k$
of a simple path $v_1v_2\cdots v_k$ and identifying $v$ with the vertex $v_{k+1}$ of a simple path
$v_{k+1}v_{k+2}\cdots v_n$, $1<k<n-1$; $G_2$ is obtained from $G_1$ by deleting $v_{k-1}v_k$ and
adding $v_1v_n$; $G_3$ is obtained from $G_1$ by deleting $v_{k+1}v_{k+2}$ and adding $v_1v_n$.
Then $Z(G_1)<Z(G_2)$ or $Z(G_1)<Z(G_3)$.
\end{lemma}

Note that the previous lemma holds also for $t = 0$ (if the vertices $u_1, u_2, \ldots, u_t$ do not exist).

\section{The second maximal Hosoya index of unicyclic graphs}

In order to get the main result, we firstly determine unicyclic graphs with the second maximal
Hosoya index in ${\cal{U}}_{n,k}$.

\begin{lemma}\label{7}
Let $U_{n,k}\in {\cal{U}}_{n,k}$ be a unicyclic graph  with the second maximal Hosoya index and
girth $k$ $(3\leq k\leq n-2)$. Then $U_{n,k}$ must be of the form $H_{i}$ $(i=1,2,3)$ (see Fig.~2),
where $s+t+k=n$.
\end{lemma}

\begin{proof}
Let $C_{k}$ be the unique cycle of $U_{n,k}$.

If there are at least three vertices on $C_{k}$ of degree greater than 2, then by Lemma \ref{5} and
Lemma \ref{6}, there exists a graph $G_1$ of the form $H_{1}$ such that $Z(U_{n,k}) < Z(G_{1}) <
Z(L_{n,k})$, which contradicts to the fact that $U_{n,k}$ has the second maximal Hosoya index.
Hence, there are at most two vertices in $C_k$ of degree at least 3 and let $u$ be a vertex of
$C_{k}$ of degree greater than~2.

If $deg (u) \geq 5$, by Lemma \ref{5} and Lemma \ref{6}, there exists a graph $G_2$ of the form
$H_2$ such that $Z(U_{n,k})< Z(G_{2})< Z(L_{n,k})$, which contradicts to the fact that $U_{n,k}$
has the second maximal Hosoya index. Hence, $3 \leq deg (u) \leq 4$. Let $v$ be a vertex of degree
greater than 2, different from $u$. We consider the following two subcases.

Assume first that $v \in V(C_{k})$. Similarly as above, we have $3 \leq deg (v) \leq 4$. If at
least one of the vertices $u$ and $v$ has degree 4, then by Lemma \ref{5}, there exists a graph
$G_3$ of the form $H_1$ such that $Z(U_{n,k})< Z(G_{3})< Z(L_{n,k})$, which is a contradiction.
Therefore, $deg (u) = deg (v) =3$ and all other vertices of $U_{n,k}$ have degree 1 or 2.

Let now $v \not \in V(C_{k})$. If there are other vertices different from $v$ and not in $C_k$ with
degree greater than 2, then by Lemma \ref{6}, we can get a graph $G_4$ of the form $H_{3}$ such
that $Z(U_{n,k})<Z(G_{4})<Z(L_{n,k})$, which is impossible. Hence, we can assume that $u$ and $v$
are the only vertices of degree greater than 2 in $U_{n,k}$. Similarly as above, we can conclude
that $deg (u) = deg (v) = 3$.
\end{proof}

\setlength{\unitlength}{0.1 in}
\begin{center}
\begin{picture}(54,11)
\put(3,8){\circle{5}}\put(1.5,8){$C_{k}$}

\put(5.3,9){\circle*{0.5}} \put(5.3,9){\line(1,0){2}}

\put(7.3,9){\circle*{0.5}}\put(8,8.6){$\cdots $}

 \put(6.7,9.5){$u_1$}\put(3.7,9){$u_0$}

\put(10.3,9){\line(1,0){2}}

\put(10.3,9){\circle*{0.5}} \put(9,9.5){$u_{t-1}$}

\put(12.3,9){\circle*{0.5}}\put(12.3,9.5){$u_t$}


\put(5.3,7){\circle*{0.5}}\put(5.3,7){\line(1,0){2}}

\put(7.3,7){\circle*{0.5}}

\put(3.7,7){$v_0$}

\put(8,6.6){$\cdots $}

 \put(6.7,6){$v_1$}

\put(10.3,7){\line(1,0){2}}

\put(10.3,7){\circle*{0.5}} \put(9,6){$v_{s-1}$}

\put(12.3,7){\circle*{0.5}}\put(12.3,6){$v_s$}

\put(2.5,2.7){$H_{1}=L^{1}_{n,k}(s,t)$}


\put(20,8){\circle{5}}\put(18.5,8){$C_{k}$}
\put(22.5,8){\circle*{0.5}}\put(23,8){$u_0$}

 \put(22.5,8){\line(1,1){2}}
\put(22.5,8){\line(1,-1){2}}

 \put(24.5,10){\circle*{0.5}}
\put(24.5,6){\circle*{0.5}}
 \put(23.5,10.6){$u_1$}
\put(24,5){$v_1$}

 \put(25.5,9.7){$\cdots$ }
\put(25,5.7){$\cdots$ }

 \put(27.5,10){\circle*{0.5}}
\put(27.5,6){\circle*{0.5}}
 \put(27.5,10){\line(1,0){2}}
\put(27.5,6){\line(1,0){2}}

 \put(29.5,10){\circle*{0.5}}
\put(29.5,6){\circle*{0.5}}

 \put(29.5,10.5){$u_{t}$}
\put(29.5,5){$v_{s}$}

 \put(26,10.6){$u_{t-1}$}
\put(26.2,5){$v_{s-1}$}

\put(21,2){$H_{2}=L^{2}_{n,k}(s,t)$}
\put(37,8){\circle{5}}\put(35,8){$C_{k}$}

\put(38.5,6){\circle*{0.5}}\put(36.7,6){$u_{0}$}

\put(38.5,6){\line(1,0){2}}\put(40,5){$u_{1}$}

\put(40.5,6){\circle*{0.5}} \put(41.3,5.7){$\cdots$}
\put(44,6){\circle*{0.5}} \put(43.5,5){$u_{l}$}

\put(49,6){\line(1,0){2}}

\put(44,6){\line(1,0){2}}\put(46,6){\circle*{0.5}}

\put(46.5,5.7){$\cdots$}

\put(49,6){\circle*{0.5}}

 \put(48.3,5){$u_{t-1}$}

 \put(51.7,5){$u_{t}$}

\put(51,6){\circle*{0.5}}

\put(44,6){\line(0,1){2}}\put(44,8){\circle*{0.5}}

\put(43.8,9){$v_1$}

\put(44,8){\line(1,0){2}}\put(46,8){\circle*{0.5}}

\put(45.8,9){$v_2$}\put(46.7,7.7){$\cdots$}

\put(49,8){\circle*{0.5}}

\put(48.3,9){$v_{s-1}$}
\put(49,8){\line(1,0){2}}\put(51,8){\circle*{0.5}}

\put(51.7,8){$v_{s}$}

\put(38,3){$H_{3}=L_{n,k}^{3}(s,t;l)$}

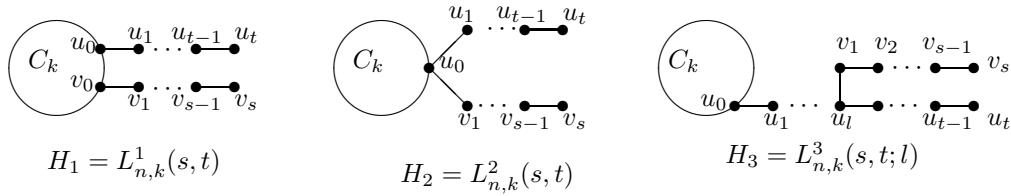
 \put(9.5,-0.8){\textbf{Figure 2.} Three classes of unicyclic
graphs.}
\end{picture}
\end{center}

Let ${\cal{L}}^{i}_{n,k}$ be the set of all unicyclic graphs  of the form like $H_{i}$ $(i=1,2,3)$
(as shown in Fig.~3),  respectively. \vspace{0.5cm}

\setlength{\unitlength}{0.1 in}
\begin{center}
\begin{picture}(54,11)
\put(3,8){\circle{5}}\put(1.5,8){$C_{k}$}

\put(5.3,9){\circle*{0.5}} \put(5.3,9){\line(1,0){2}}

\put(7.3,9){\circle*{0.5}}\put(8,8.6){$\cdots $}

 \put(6.7,9.5){$u_1$}\put(3.7,9){$u_0$}

\put(10.3,9){\line(1,0){2}}

\put(10.3,9){\circle*{0.5}}

\put(12.3,9){\circle*{0.5}}\put(12.3,9.5){$u_{n-k-2}$}


\put(5.3,7){\circle*{0.5}}\put(5.3,7){\line(1,0){4}}

\put(7.3,7){\circle*{0.5}}\put(9.3,7){\circle*{0.5}}

\put(6.7,6){$y$} \put(3.7,7){$v_0$}

\put(8.7,6){$x$}

%
%
%
%

\put(2.5,2.7){$L^{1}_{n,k}$}


\put(20,8){\circle{5}}\put(18.5,8){$C_{k}$}
\put(22.5,8){\circle*{0.5}}\put(23,8){$u_0$}

 \put(22.5,8){\line(1,1){2}}
\put(22.5,8){\line(1,-1){2}}

 \put(24.5,10){\circle*{0.5}}
\put(24.5,6){\circle*{0.5}}
 \put(23.5,10.6){$u_1$}
\put(24,5){$y$}

 \put(25.5,9.7){$\cdots$ }

 \put(27.5,10){\circle*{0.5}}
 \put(27.5,10){\line(1,0){2}}
\put(24.5,6){\line(1,0){2}}

 \put(29.5,10){\circle*{0.5}}
\put(26.5,6){\circle*{0.5}}

 \put(29.5,10.5){$u_{n-k-2}$}
\put(26.5,5){$x$}


\put(21,2){$L^{2}_{n,k}$}
\put(37,8){\circle{5}}\put(35,8){$C_{k}$}

\put(38.5,6){\circle*{0.5}}\put(36.7,6){$u_{0}$}

\put(38.5,6){\line(1,0){7}}\put(40,5){$u_{1}$}

\put(40.5,6){\circle*{0.5}}
\put(43,6){\circle*{0.5}} \put(42.5,5){$u_{2}$}

\put(49,6){\line(1,0){2}}

\put(44,6){\line(1,0){2}}\put(46,6){\circle*{0.5}}

\put(46.5,5.7){$\cdots$}

\put(49,6){\circle*{0.5}}


 \put(51.1,5){$u_{n-k-2}$}

\put(51,6){\circle*{0.5}}

\put(40.5,6){\line(1,1){2}}

\put(42.5,8){\circle*{0.5}}

\put(42.5,8){\line(1,0){2}} \put(44.5,8){\circle*{0.5}}
\put(44.5,8.5){$y$}\put(42.5,8.5){$x$}

%
%
%
%
%
%
%

\put(38,3){$L^{3}_{n,k}$}

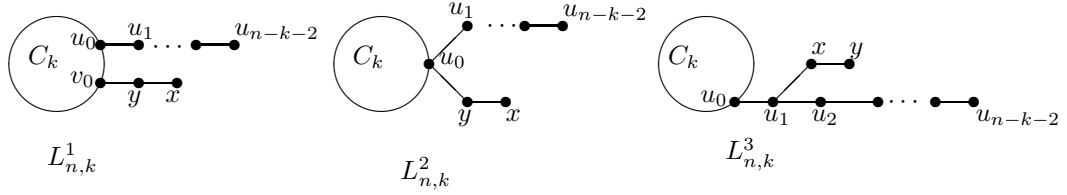
 \put(9.5,-0.8){\textbf{Figure 3.}  The extremal unicyclic
graphs.}
\end{picture}
\end{center}

\begin{lemma}\label{8} Let $n \geq 10$.
\begin{enumerate}[($i$)]
\item If $k>3$, $L^{i}_{n,k}$ $(i=1,2,3)$ (see Fig.~3) is the unique graph with the maximal
Hosoya index in ${\cal{L}}^{i}_{n,k}$ $(i=1,2,3)$, respectively. In $L^{1}_{n,k}$, the vertices
$u_0$ and $v_0$ are adjacent;

\item If $k=3$, $L^{3}_{n,3}(2,n-5;3)$ (see Fig. ~2) and $L_{n,3}^{3}$ (see Fig.~3) are the graphs with the maximal
Hosoya index in ${\cal{L}}^{3}_{n,3}$; $L^{i}_{n,3}$ $(i=1,2)$ (see Fig.~3) is the unique graph
with the maximal Hosoya index in ${\cal{L}}^{i}_{n,3}$ $(i=1,2)$, respectively. In $L^{1}_{n,3}$,
the vertices $u_0$ and $v_0$ are adjacent.
\end{enumerate}
\end{lemma}

\begin{proof} We distinguish the following three cases.

{\bf Case 1.} For graphs in ${\cal{L}}^{1}_{n,k}$.

By formula (\ref{formula4}), we have
\begin{eqnarray*}
Z(L^{1}_{n,k}(s,t))&=&Z(L^{1}_{n,k}(s,t)-v_{0}v_{1})+Z(L^{1}_{n,k}(s,t)-v_{0}-v_{1})\\
&=&Z(P_{s})Z(L_{n-s,k})+Z(P_{s-1}) Z(T_{n-s-1})\\
&\leq&Z(P_{s})Z(L_{n-s,k})+Z(P_{s-1}) Z(P_{n-s-1}),
\end{eqnarray*}
where $T_{n-s-1}$ is a tree of order $n-s-1$. The equality holds if
and only if $T_{n-s-1}\cong P_{n-s-1}$ and consequently $u_0$ and
$v_0$ must be adjacent in $L^{1}_{n,k}(s,t)$ (see Fig.~2). In the
following, we assume $u_0$ and $v_0$ are adjacent in
$L^{1}_{n,k}(s,t)$.
\begin{eqnarray*}
Z(L^{1}_{n,k}(s,t))&=&Z(P_{s}) Z(L_{n-s,k})+Z(P_{s-1}) Z(P_{n-s-1})\\
&=&F_{s+1}(F_{n-s+1}+F_{k-1}F_{n-s-k+1})+F_{s}F_{n-s}\\
&=&F_{s+1}F_{n-s+1}+F_{s}F_{n-s}+F_{k-1}\cdot F_{s+1}F_{n-s-k+1}\\
&=&F_{n+1} + F_{k-1}\cdot F_{s+1}F_{n-s-k+1}.
\end{eqnarray*}

If $s$ is odd, $Z(L^{1}_{n,k}(s,t))$ is strictly increasing with respect to $s$. If $s$ is even,
$Z(L^{1}_{n,k}(s,t))$ is strictly decreasing with respect to $s$. Therefore, it follows that
$$
Z(L^{1}_{n,k}(2,n-k-2))>Z(L^{1}_{n,k}(4,n-k-4))>\cdots>Z(L^{1}_{n,k}(2m,2m+i))
>
$$
$$
> Z(L^{1}_{n,k}(2m-1+2l,2m+i+1-2l))>\cdots>Z(L^{1}_{n,k}(3,n-k-3))>Z(L^{1}_{n,k}(1,n-k-1)),
$$
where $n - k = 4m + i$, $i \in \{0, 1, 2, 3\}$ and $l = \lfloor \frac{i}{2} \rfloor$. Finally, we
get that $L_{n,k}^1$ is the only graph with the maximal Hosoya index in ${\cal{L}}^{1}_{n,k}$.

{\bf Case 2.} For graphs in ${\cal{L}}^{2}_{n,k}$.

We can easily get the result from Lemma \ref{5}.
It follows that $L_{n,k}^2$ is the only graph with the maximal Hosoya index in ${\cal{L}}^{2}_{n,k}$.

{\bf Case 3.} For graphs in ${\cal{L}}^{3}_{n,k}$.

In order to prove the $L_{n,k}^3$ is the only graph with the maximal Hosoya index in
${\cal{L}}^{3}_{n,k}$, it suffices to verify that
$Z(L_{n,k}^{3}(s,t;l))\leq Z(L^{3}_{n,k})$.
For $s\geq 1$, by Lemma \ref{5}, there is an integer $l_0$ such that $Z(L_{n,k}^{3}(2,n-k-2;l_0))\geq
Z(L_{n,k}^{3}(s,t;l))$. For convenience, we denote
$G=L_{n,k}^{3}(2,n-k-2;l_0)$.

By formula (\ref{formula6}), we have
\begin{eqnarray*}
Z(G)&=&Z(G-v_{2})+Z(G-v_{2}-v_{1})\\
&=&Z(G-v_{2}-v_{1})+Z(G-v_{2}-v_{1}-u_{l_0})+Z(G-v_{2}-v_{1})\\
&=&2Z(L_{n-2,k})+Z(P_{l_0-1}) Z(L_{n-l_0-2,k}).
\end{eqnarray*}
\begin{eqnarray*}
Z(L^{3}_{n,k})&=&Z(L^{3}_{n,k}-x)+Z(L^{3}_{n,k}-x-y)\\
&=&Z(L^{3}_{n,k}-x-y)+Z(L^{3}_{n,k}-x-y-u_{1})+Z(L^{3}_{n,k}-x-y)\\
&=&2Z(L_{n-2,k})+Z(C_{k}) Z(P_{n-k-3}).
\end{eqnarray*}
Therefore, it follows that
\begin{eqnarray*}
Z(L^{3}_{n,k})-Z(G)&=&Z(C_{k}) Z(P_{n-k-3})-Z(P_{l_0-1})
Z(L_{n-l_0-2,k})\\
&=&Z(C_{k}) Z(P_{n-k-3})-Z(P_{l_0-1}) [Z(C_{k})Z(P_{n-k-l_0-2})+Z(P_{k-1})Z(P_{n-k-l_0-3})]\\
&=&Z(C_{k}) F_{n-k-2}-F_{l_0}[Z(C_{k})F_{n-k-l_0-1}+F_{k}F_{n-k-l_0-2}]\\
&=&Z(C_{k})[F_{n-k-2}-F_{l_0}F_{n-k-l_0-1}]-F_{l_0}F_{k}F_{n-k-l_0-2}\\
&=&Z(C_{k})F_{l_0-1}F_{n-k-l_0-2}-F_{l_0}F_{k}F_{n-k-l_0-2}\\
&=&[F_{k+1}F_{l_0-1}+F_{k-1}F_{l_0-1}-F_{l_0}F_{k}]F_{n-k-l_0-2}\\
&=&[2F_{k-1}F_{l_0-1}-F_{k}F_{l_0-2}]F_{n-k-l_0-2}\\
&\geq&[2F_{k-1}F_{l_0-1}-F_{k}F_{l_0-1}]F_{n-k-l_0-2}\\
&=&[F_{k-1}F_{l_0-1}-F_{k-2}F_{l_0-1}]F_{n-k-l_0-2}\\
&=&F_{k-3}F_{l_0-1}F_{n-k-l_0-2}\geq 0.
\end{eqnarray*}
As above, it is easy to find that the equality holds if and only if $s = 2$ and $l_0 = n - k - 2$, or $k=3$, $s=2$ and $l_0=3$.

For $k = 3$, we have
$$Z(L^3_{n,3})=2Z(L_{n-3,3})+4F_{n-5};$$
$$Z(L^3_{n,3}(2,n-5;3))=2Z(L_{n-2,3})+2Z(L_{n-5,3}).$$
It is easy to verify that
\begin{eqnarray*}
Z(L^3_{n,3})-Z(L^3_{n,3}(2,n-5;3))&=&4F_{n-5}-2Z(L_{n-5,3})\\
&=&4F_{n-5}-2F_{n-4}-2F_{n-7}\\
&=&0.
\end{eqnarray*}
This completes the proof.
\end{proof}

\setlength{\unitlength}{0.1 in}
\begin{center}
\begin{picture}(30,11)
\put(7,8){\circle{5}}\put(4.6,8){$C_{n-2}$}
\put(11.3,9){\circle*{0.5}} \put(9.3,9){\circle*{0.5}}
\put(7.7,9){$v_0$}

 \put(9.3,9){\line(1,0){2}}

\put(11.3,7){\circle*{0.5}} \put(9.3,7){\circle*{0.5}}
\put(7.7,7){$u_0$}\put(13.3,7){\circle*{0.5}}

\put(9.3,7){\line(1,0){4}}

\put(6.5,3){$U'_{n}$}

\put(21,8){\circle{5}}\put(18.6,8){$C_{n-3}$}

\put(23.2,9){\circle*{0.5}} \put(23.2,7){\circle*{0.5}}
\put(21.3,8.8){$v_0$} \put(21.3,6.8){$u_0$}
\put(23.2,9){\line(1,0){2}}\put(25.2,9){\circle*{0.5}}

\put(23.2,7){\line(1,0){4}}\put(25.2,7){\circle*{0.5}}\put(27.2,7){\circle*{0.5}}
\put(21.5,3){$U''_{n}$}

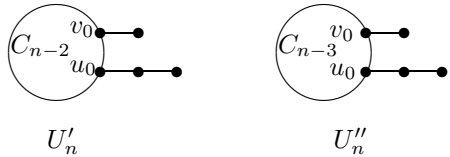
 \put(-3,-0.1){\textbf{Figure 4.} Extremal  graphs for $k = n - 2$ and $k = n - 3$.}
\end{picture}
\end{center}

\begin{lemma} \label{9}
Let $U_{n,k}\in {\cal{U}}_{n,k}\backslash \{L_{n,k} \}$ be a unicyclic graph with girth $k$ ($3\leq
k\leq n-2$) and $n \geq 10$.
\begin{enumerate}[($i$)]
\item If $k=n-2$, then $Z(U_{n,k}) \leq Z(U'_{n}) = F_{n+1} + F_{n-3}$,
with equality if and only if $U_{n,k}\cong Z(U'_{n})$.
\item If $k=n-3$, then $Z(U_{n,k}) \leq Z(U''_{n}) = F_{n+1} + 2F_{n-4}$,
with equality if and only if $U_{n,k}\cong Z(U''_{n})$.
\item If $3\leq k\leq \frac{n-1}{2}$ and $k\neq\frac{n-2}{2}$, then
$Z(U_{n,k})\leq F_{n+1}+2F_{k-1} F_{n-k-1}$ for $k$ odd, with equality if and only if $U_{n,k}\cong
L_{n,k}^{1}$; $Z(U_{n,k})\leq 2F_{n-1}+F_{k-1} F_{n-k+1}+F_{k+1} F_{n-k-2}$ for $k$ even, with
equality if and only if $U_{n,k}\cong L_{n,k}^{3}$.
\item If $k=\frac{n-2}{2}$, then
$Z(U_{n,k})\leq F_{n+1}+2F_{k-1}F_{n-k-1}$, with equality if and only $U_{n,k}\cong L_{n,k}^{1}$ or
$U_{n,k}\cong L_{n,k}^{3}$.
\item If $\frac{n-1}{2}<k\leq n-4$, then $Z(U_{n,k})\leq
F_{n+1}+2F_{k-1} F_{n-k-1}$ for $k$ and $n$ having the same parity, with equality if and only if
$U_{n,k}\cong L_{n,k}^{1}$; otherwise, $Z(U_{n,k})\leq 2F_{n-1}+F_{k-1}
F_{n-k+1}+F_{k+1}F_{n-k-2}$, with equality if and only if $U_{n,k}\cong L_{n,k}^{3}$.
\end{enumerate}
The graphs $L_{n,k}^{1}, L_{n,k}^{3}, U'_n, U''_n$ are shown in
Fig.~3 and Fig.~4 and $u_0$ is adjacent to $v_0$ in $U'_n, U''_n$.
\end{lemma}

\begin{proof}
For $k=n-2$, from Lemma \ref{8} the unique extremal graph is $U'_{n} \cong L^{1}_{n,n-2}$. For
$k=n-3$, from Lemma \ref{8} and direct verification we have that the unique extremal graph is
$U''_{n} \cong Z(L^1_{n,n - 3})$. For $k=3$, we can easily verify that $Z(L^1_{n,3})>Z(L^2_{n,3})$,
and $Z(L^1_{n,3})>Z(L^3_{n,3}(2,n-5;3))$.

For $3 < k \leq n-4$, by Lemma \ref{8}, it suffices to
compare the values $Z(L^{1}_{n,k})$, $Z(L^{2}_{n,k})$ and
$Z(L^{3}_{n,k})$.

By formula (\ref{formula6}), we have
\begin{eqnarray*}
Z(L^{1}_{n,k})&=&2Z(L_{n-2,k})+Z(P_{n-3}),\\
Z(L^{2}_{n,k})&=&2Z(L_{n-2,k})+Z(P_{k-1})Z(P_{n-k-2}),\\
Z(L^{3}_{n,k})&=&2Z(L_{n-2,k})+Z(C_{k})Z(P_{n-k-3}).
\end{eqnarray*}

It is not difficult to verify that $Z(L^{1}_{n,k}) - Z(L^{2}_{n,k}) = F_{n-2} - F_k F_{n-k-1} >
F_{n-2} - F_3 F_{n-4} = F_{n-5} > 0$ for $n > 5$.
\begin{eqnarray*}
Z(L^{1}_{n,k})-Z(L^{3}_{n,k})&=&Z(P_{n-3})-Z(C_{k})Z(P_{n-k-3})\\
&=&F_{n-2}-(F_{k+1}+F_{k-1})F_{n-k-2}\\
&=&F_{k}F_{n-k-1}+F_{k-1}F_{n-k-2}-F_{k+1}F_{n-k-2}-F_{k-1}F_{n-k-2}\\
&=&F_{k}F_{n-k-1}-F_{k+1}F_{n-k-2}\\
&=&(-1)^{k+1}F_{n-2k-2} \qquad \qquad \mbox {by formula \eqref{fi3}}.
\end{eqnarray*}

If $3< k\leq \frac{n-1}{2}$ and $k\neq \frac{n-2}{2}$, then $n \geq
2k + 1$. For $n > 2k + 2$ it follows that $F_{n-2k-2}>0$, and for $n
= 2k + 1$ it follows that $F_{-1} = 1 > 0$. Therefore, in both cases
if $k$ is odd, then $Z(L^{1}_{n,k})-Z(L^{3}_{n,k})>0$; if $k$ is
even, then $Z(L^{1}_{n,k})-Z(L^{3}_{n,k})<0$.

If $k=\frac{n-2}{2}$, then $F_{n-2k-2}=F_0=0$ and $Z(L^{1}_{n,k})=Z(L^{3}_{n,k})$.

If $\frac{n-1}{2}<k\leq n-4$, then we have
$$(-1)^{k+1} F_{n-2k-2} = (-1)^{k+1} (-1)^{2k+3-n}F_{2k+2-n} = (-1)^{n+k} F_{2k+2-n}.$$
So, if $k$ and $n$ have the same parity, then
$Z(L^{1}_{n,k})-Z(L^{3}_{n,k})>0$; otherwise, $Z(L^{1}_{n,k})-Z(L^{3}_{n,k})<0$.

This completes the proof.
\end{proof}

\section{Ordering of unicyclic graphs with respect to Hosoya index}

In this section, we extend the result in \cite{Ou jianping} and order the first $n - 1$ unicyclic
graphs with respect to the Hosoya index for $n \geq 11$. For $5 \leq n \leq 10$, using exhaustive
computer search among all unicyclic graphs (with the help of Nauty~\cite{Nauty}), we have the
following:

If $n=5$, then $Z(C_5) = 11 > Z(L_{5,4}) = Z(L_{5,3}) = 10 > Z(U'_5) = 9.$

If $n=6$, then $Z(C_6) = 18 > Z(L_{6,4}) = 17 > Z(L_{6,5}) = Z(L_{6,3}) = 16 > Z(U''_6) = Z(U'_6) =
15.$

If $n=7$, then $Z(C_7) = 29 > Z(L_{7,4})=Z(L_{7,5}) = 27 > Z(L_{7,3})=Z(L_{7,6}) = 26 >
Z(L^1_{7,3}) = Z (U''_7) = 25.$

If $n=8$, then $Z(C_8) = 47>Z(L_{8,4})=Z(L_{8,6}) = 44
>Z(L_{8,5}) = 43 >Z(L_{8,3})=Z(L_{8,7})=Z(L^1_{8,4}) = 42.$

If $n=9$, then $Z(C_9) = 76>Z(L_{9,4})=Z(L_{9,7}) = 71 >Z(L_{9,6})=Z(L_{9,5}) = 70
>Z(L_{9,8})=Z(L_{9,3})=Z(L^3_{9,4}) = 68.$

If $n=10$, then $Z(C_{10}) = 123 > Z(L_{10,4})=Z(L_{10,8}) = 115 > Z(L_{10,6}) = 114
>Z(L_{10,7})=Z(L_{10,5}) = 113 > Z(L_{10,9}) = Z(L_{10,3}) = 110 > Z (L^{3}_{10, 4}) = Z (L^{1}_{10, 4})
= Z(L^{1}_{10, 6}) = 109.$

All other unicyclic graphs have strictly smaller Hosoya index.

\begin{theorem}
For $n \geq 11$, the first $n-1$ largest Hosoya indices of unicyclic graphs are: $Z(C_n) >
Z(L_{n,4}) > Z(L_{n,6}) > \cdots > Z(L_{n, 2m}) > Z(L_{n, 2m+1+2l})
> \cdots > Z(L_{n,5}) > Z(L_{n,3}) > Z(L^{3}_{n, 4})$,
where $n = 4m + i$, $i \in \{0, 1, 2, 3\}$, $l = \lfloor \frac{i}{2} \rfloor$, and $Z (L_{n,k}) = Z
(L_{n, n-k+2})$.
\end{theorem}

\begin{proof}
Let ${\cal{U}}_{n}$ be the set of all connected unicyclic graphs on $n$ vertices. Since
${\cal{U}}_{n}$ is the union of all ${\cal{U}}_{n,k}$, where $k = 3, 4, \ldots, n$, we need to
order the extremal graphs $L^{1}_{n,k}$, $L^{3}_{n,k}$ and $L_{n,k}$, $k = 3, 4, \ldots, n$ based
on the Hosoya index. It follows that for $n \geq 11$, the second to the $(n-2)$-th largest Hosoya
indices are exactly graphs $L_{n, k}$, $k = 3, 4, \ldots, n - 1$, while the $(n-1)$-th largest
Hosoya index is achieved uniquely by the unicyclic graph~$L^{3}_{n, 4}$. \vspace{0.2cm}

\textbf{Claim 1.} $ Z(L^{1}_{n,4}) > Z(L^{1}_{n,6}) > \cdots > Z(L^{1}_{n, 2m}) > Z(L^{1}_{n, 2m-1
+ 2l}) > \cdots > Z(L^{1}_{n,5}) > Z(L^{1}_{n,3}),$ where $n = 4m + i$, $i \in \{0, 1, 2, 3\}$, and
$l = \lfloor \frac{i}{2} \rfloor$.

First note that $Z(L^{1}_{n,k}) = Z(L^{1}_{n,n-k})$, so we can assume $k \leq \frac n2$. From the
proof of Lemma~\ref{9} and Lemma \ref{4}, we have
\begin{eqnarray*}
Z(L^{1}_{n,k}) - Z(L^{1}_{n,k-1}) &=& 2 \cdot \left( Z (L_{n-2,k}) - Z (L_{n-2,k-1}) \right )\\
&=& 2 \cdot \left( F_{k-1} F_{n-k-1} - F_{k-2}F_{n-k} \right ) \\
&=& 2 (-1)^{k} \cdot F_{n-2k+1} \quad \mbox {by formula \eqref{fi3}}.
\end{eqnarray*}
Hence, if $k$ is even, $Z(L^{1}_{n,k}) - Z(L^{1}_{n,k-1})>0$; if $k$
is odd, $Z(L^{1}_{n,k}) - Z(L^{1}_{n,k-1})<0$.
\begin{eqnarray*}
Z(L^{1}_{n,k}) - Z(L^{1}_{n,k-2}) &=& Z(L^{1}_{n,k}) -
Z(L^{1}_{n,k-1})+Z(L^{1}_{n,k-1}) - Z(L^{1}_{n,k-2})\\
&=& 2 (-1)^{k} \cdot \left( F_{n-2k+1} -F_{n-2k+3} \right ).
\end{eqnarray*}
Hence, if $k$ is even, $Z(L^{1}_{n,k}) - Z(L^{1}_{n,k-2})<0$; if $k$
is odd, $Z(L^{1}_{n,k}) - Z(L^{1}_{n,k-2})>0$.\vspace{0.2cm}

\textbf{Claim 2. }$ Z(L_{n,4}) > Z(L_{n,6}) > \cdots > Z(L_{n, 2m}) > Z(L_{n, 2m+1+2l}) > \cdots >
Z(L_{n,5}) > Z(L_{n,3}), $ where $n = 4m + i$, $i \in \{0, 1, 2, 3\}$, and $l = \lfloor \frac{i}{2}
\rfloor$.

For the extremal graphs $L_{n,k}$, we have $Z (L_{n, k}) = Z (L_{n, n - k + 2})$ and
\begin{eqnarray*}
Z(L_{n,k}) - Z(L_{n,k-1}) &=& F_{k-1} F_{n-k+1} -F_{k-2}F_{n-k+2}\\
&=& (-1)^{k} \cdot F_{n-2k+3}.
\end{eqnarray*}

Using the formula (\ref{fib}), we complete the result.

\vspace{0.2cm} \textbf{Claim 3. } $Z(L^{3}_{n, k}) \leq Z(L^{3}_{n, 4})$
for $3\leq k\leq n-4$, with equality if and only if $k = 4$.

Note that $Z(L^{3}_{n, k}) = 2F_{n-1}+F_{k-1} F_{n-k+1}+F_{k+1}F_{n-k-2}$ and $Z (L^{3}_{n, 4}) =
2F_{n-1} + F_3 F_{n - 3} + F_5 F_{n-6}$. From formula (\ref{fib}), we have
$$
F_3 F_{n-3} > F_{k-1} F_{n-(k-1)} \qquad \mbox{ for } \ 3 \leq k \leq n - 4, \ \ k \neq 4,
$$
and
$$
F_5 F_{n - 1 - 5} > F_{k+1} F_{n - 1 - (k + 1)} \qquad \mbox{ for } \ 3 \leq k \leq n - 4, \ \ k \neq 4, n - 5.
$$
For the special case $k = n - 5$, we have $Z (L^{3}_{n,n-5}) = 2 F_{n-1}+F_{n-6} F_6 + F_{n-4} F_3$
and obviously $Z (L^{3}_{n, 4}) - Z (L^{3}_{n, n-5}) = F_3 F_{n-5} - F_{4} F_{n-6} = 2 F_{n-5} - 3
F_{n-6} = F_{n-9} > 0$. Therefore, it follows that $Z(L^{3}_{n, k}) < Z(L^{3}_{n, 4})$, for $3\leq
k\leq n-4$ and $k \neq 4$.


\vspace{0.2cm}\textbf{ Claim 4.} The cycle $C_n$ has the largest Hosoya index among unicyclic graphs \cite{Ou jianping}. By using
the above three claims, in order to prove that the extremal graphs $L_{n, k}$ $(k = 3, 4, \ldots, n
- 1)$ are the next $n - 3$ unicyclic graphs with the largest Hosoya index, we need to compare the
following Hosoya indices
\begin{eqnarray*}
Z (L_{n, k}) = Z (L_{n, n - k + 2}) &=& F_{n+1} + F_{k-1} F_{n - k + 1},\\
Z (L^{1}_{n, 4}) = Z (L^{1}_{n, n - 4}) &=& F_{n+1} + 4 F_{n - 5}, \\
Z (L^{3}_{n, 4}) &=& 2F_{n-1} + 2 F_{n - 3} + 5 F_{n-6},\\
Z (U'_n) &=& F_{n+1} + F_{n-3} ,\\
Z (U''_n) &=& F_{n+1} + 2F_{n-4}.
\end{eqnarray*}

It is easy to verify that $Z (L^{1}_{n, 4}) > Z (U''_n) > Z (U'_n)$ and $Z (L^{3}_{n, 4}) > Z
(L^{1}_{n, 4})$. Therefore, we need to estimate the difference
$$
\Delta = Z (L_{n, k}) - Z (L^{3}_{n, 4}) = F_{n + 1} + F_{k - 1} F_{n - k + 1} - (2F_{n-1} + 2 F_{n - 3} + 5 F_{n-6}).
$$

Note that for $k \geq 3$, $F_{k - 1} F_{n - (k - 1)} \geq F_{2} F_{n
- 2} = F_{n - 2}$. Hence,
\begin{eqnarray*}
\Delta &\geq& F_{n+1} + F_{n-2} - 2F_{n-1} - 2 F_{n - 3} - 5 F_{n-6} \\
&=& 2 F_{n-2} - 2 F_{n-3} - 5 F_{n-6} \\
&=& 2 F_{n-5} - 3 F_{n-6} \\
&=& 2 F_{n-7} - F_{n-6} \\
&=& F_{n-9} > 0.
\end{eqnarray*}

This completes the proof.
%
%
%
%
%
%
\end{proof}


\end{document}